\documentclass[12pt]{article}
\usepackage{amssymb}
\usepackage{amscd}
\newtheorem{theorem}{Theorem}[section]
\newtheorem{proposition}[theorem]{Proposition}
\newtheorem{lemma}[theorem]{Lemma}
\newtheorem{remark}[theorem]{Remark}
\newtheorem{corollary}[theorem]{Corollary}

\newtheorem{list1}[theorem]{List}
\newtheorem{list2}[theorem]{List}
\newcommand{\classno}[1]{\medskip 2000 Mathematics Subject Classification~~ #1}

\def\carre {{\vrule height5pt width5pt depth0pt}}
\def\qed{{\hfill\carre\vskip1em}}
\newenvironment{proof}{\medskip\goodbreak\noindent{\it Proof.~}}{\qed}
\newenvironment{acknowledgements}{\medskip\goodbreak{\it Acknowledgements.~}}{}
\newcommand{\affiliationone}[1]{\noindent#1}
\def\Hom{\mathrm{Hom}}
\def\ad{\mathrm{ad}}
\title
{On the irreducibility of the commuting variety of a symmetric pair
associated to a parabolic subalgebra with abelian unipotent radical}
\author{Herv\'e Sabourin and Rupert W.T. Yu}

\begin{document}
\maketitle

\begin{abstract}
In this paper, we study the commuting variety of
symmetric pairs associated to parabolic subalgebras with abelian
unipotent radical in a simple complex Lie algebra. By using the ``cascade''
construction of Kostant, we construct a Cartan subspace which in turn provides,
in certain cases, useful information on the centralizers of 
non $\mathfrak{p}$-regular semisimple elements. 
In the case of the rank $2$ symmetric pair 
$(\mathrm{so}_{p+2},\mathrm{so}_{p}\times \mathrm{so}_{2})$,
$p\geq 2$, this allows us to apply induction, in view of previous results of 
the authors~\cite{SY}, and reduce the problem of the irreducibility of 
the commuting variety to the consideration of evenness of $\mathfrak{p}$-distinguished 
elements. Finally, via the correspondence of Kostant-Sekiguchi, we check
that in this case, $\mathfrak{p}$-distinguished elements are indeed
even, and consequently, the commuting variety is irreducible.
\end{abstract}

\classno{17B20, 14L30}

\section{Introduction and notations}\label{intro}

Let $\mathfrak{g}$ be a complex simple Lie algebra
and $\theta$ an involutive automorphism of $\mathfrak{g}$. 
Let $\mathfrak{g} = \mathfrak{k} \oplus
\mathfrak{p}$ be the decomposition of $\mathfrak{g}$ into 
eigenspaces with respect to $\theta$, where 
$\mathfrak{k} = \{ X \in \mathfrak{g} \ | \ \theta(X) = X \}$, 
$\mathfrak{p} = \{ X \in \mathfrak{g} \ | \ \theta(X) = -X \}$. 
In this case, we say that $(\mathfrak{g},\mathfrak{k})$ is a {\it symmetric pair}.

Let $G$ be the adjoint group of $\mathfrak{g}$
and $K$ the connected algebraic subgroup of $G$ whose Lie algebra is
$\mathfrak{k}$. 

Let $\mathfrak{a}$ be a maximal abelian subspace of $\mathfrak{p}$ consisting
of semisimple elements. Any such subspace is called a {\it Cartan subspace} of 
$\mathfrak{p}$. All the Cartan subspaces are $K$-conjugate. Its dimension is
called the {\it rank} of the symmetric pair $(\mathfrak{g}, \mathfrak{k})$. 

We define the commuting variety of  
$(\mathfrak{g},\mathfrak{k})$ as the following set :
$$
C(\mathfrak{p}) = \{ (x,y) \in \mathfrak{p} \times \mathfrak{p} \  | \ [x,y] = 0
\}.
$$

We may also consider the commuting variety 
$C(\mathfrak{g})$ of $\mathfrak{g}$, defined in the same
way. Richardson proved in \cite{R} that, if $\mathfrak{h}$ is a Cartan
subalgebra of $\mathfrak{g}$, then 
$C(\mathfrak{g}) = \overline{G.(\mathfrak{h} \times \mathfrak{h})}$. 
In particular, the commuting variety $C(\mathfrak{g})$ is
an irreducible algebraic variety.

On the other hand, the commuting variety of any semisimple
symmetric pair is not irreducible in general. Panyushev showed
in \cite{P1}
that in the case of the symmetric pair $(\mathrm{sl}_{n},
\mathrm{gl}_{n-1})$, $n > 2$, associated to the
involutive automorphism, defined via
conjugation by the diagonal matrix 
$\mathrm{diag}(-1, \dots,-1,1)$. The corresponding commuting variety has 
three irreducible components of dimension, respectively, $2n-1$, $2n-2$, 
$2n-2$.

Nevertheless, in some cases, the irreducibility problem has been
solved.
\begin{itemize}
\item[-] 
As an obvious consequence of the classical case proved by Richardson,
the symmetric pair $(\mathfrak{g} \times \mathfrak{g}, \Delta(\mathfrak{g}))$, 
associated to the automorphism $(X,Y)\mapsto (Y,X)$, has an irreducible
commuting variety. 
\item[-] If the rank of the symmetric pair 
$(\mathfrak{g},\mathfrak{k})$ is equal to the 
semisimple rank of $\mathfrak{g}$ (called the maximal rank case), then 
Panyushev proved in \cite{P1} that the corresponding commuting variety is irreducible.
\item[-] The rank $1$ case has been considered independently by the authors \cite{SY} 
and Panyushev \cite{P2}. In this case, it has been proved that
$(\mathrm{so}_{m+1}, \mathrm{so}_{m} )$ is the only symmetric pair 
whose commuting variety is irreducible.
\item[-] In \cite{P2}, Panyushev proves the  
irreducibility of the commuting variety for the symmetric pairs  
$(\mathrm{sl}_{2n}, \mathrm{sp}_{2n})$ and $(E_6,F_4)$.
\end{itemize}

For a symmetric pair of rank strictly larger than one, we observe that
due to the rank $1$ case, the inductive arguments used by Richardson in 
the classical case \cite{R} do not apply. However, if $\mathfrak{a}$
is a Cartan subspace, then it is well-known that 
$C_{0} = \overline{K.(\mathfrak{a}\times \mathfrak{a})}$ is the unique irreducible 
component of $C(\mathfrak{p})$ of maximal dimension, which is 
equal to $\dim \mathfrak{p} + \dim \mathfrak{a}$. The main problem is 
therefore to determine if there exist components other than the maximal one.

In \cite{P2}, it has been conjectured that $C(\mathfrak{p})$ is irreducible
if the rank of the symmetric pair is greater than or equal to $2$.

In this paper, we consider 
symmetric pairs $(\mathfrak{g}, \mathfrak{k})$ 
associated to parabolic subalgebras with abelian
unipotent radical. In this case, by using the ``cascade''
construction of Kostant, we obtain a particular Cartan subspace of 
$\mathfrak{p}$ which, in certain cases, turns out to be a useful tool
for the description of symmetric
subpairs of $(\mathfrak{g}, \mathfrak{k})$, associated to
centralizers of semisimple elements of $\mathfrak{p}$.   

The information on the centralizers of semisimple 
elements of $\mathfrak{p}$ allows us to apply the
induction step of Richardson in the case of
the rank $2$ symmetric pair 
$(\mathrm{so}_{p+2}, \mathrm{so}_{p} \times \mathrm{so}_{2})$, 
$p \geq 2$. As a result, we obtain that all the irreducible components 
of $C(\mathfrak{p})$ other than $C_{0}$ are necessarily related to
$K$-orbits of $\mathfrak{p}$-distinguished elements in $\mathfrak{p}$. 

By using the Kostant-Sekiguchi correspondence, we are
able to establish that in the case of symmetric pair 
$(\mathrm{so}_{p+2}, \mathrm{so}_{p} \times \mathrm{so}_{2})$,
$p\geq 2$, every $\mathfrak{p}$-distinguished elements is even.
We prove finally that this condition is sufficient for
the irreducibility of the commuting variety
of $(\mathrm{so}_{p+2}, \mathrm{so}_{p} \times \mathrm{so}_{2})$,
$p\geq 2$.

We shall conserve the notations above in the sequel.
The reader may refer to \cite{T} for basic definitions
and properties of symmetric pairs.
The paper is organized as follows.
In Sections \ref{parabolic} and \ref{centralizers}, 
we study symmetric pairs associated to parabolic subalgebras with abelian
unipotent radical and the description of subpairs
associated to centralizers of semisimple elements.
We recall some well-known results relating sheets
and the commuting variety in Section~\ref{commuting}.
In particular, we prove that an even nilpotent element
of $\mathfrak{p}$ belongs to a sheet containing non-zero semisimple
elements. Section~\ref{ranktwo} is dedicated to the rank $2$ symmetric pair 
$(\mathrm{so}_{p+2}, \mathrm{so}_{p} \times \mathrm{so}_{2})$, 
$p \geq 2$. Finally, Section~\ref{computations} describes
certain computations used in Section~\ref{centralizers}.

\section{Symmetric pair associated to a parabolic subalgebra with 
abelian unipotent radical}\label{parabolic}

Let us fix a Cartan subalgebra $\mathfrak{h}$ of $\mathfrak{g}$
and a Borel subalgebra $\mathfrak{b}$ containing $\mathfrak{h}$.
Denote by $R\supset R^{+}\supset \Pi$ the corresponding 
set of roots, positive roots and simple roots. Let us also
fix root vectors $X_{\alpha}$, $\alpha\in R$, and for $\alpha\in R$, 
we set $\mathfrak{g}_{\alpha}=\mathbb{C} X_{\alpha}$.

For any subset $S\subset \Pi$, we denote
by $R_{S}=R\cap \mathbb{Z}S$. Then 
$R_{S}$ is a root system in the vector subspace that
it spans in $\mathfrak{h}^{*}$. Moreover $S$ is a base
of simple roots for $R_{S}$. 

Set $R_{S}^{\pm} = R_{S}\cap R^{\pm}$, and
$$
\mathfrak{q}_{S}=\mathfrak{h}\oplus \bigoplus_{\alpha\in R_{S}\cup R^{+}} \mathfrak{g}_{\alpha}
$$
the standard parabolic subalgebra associated to $S$. 
The unipotent radical $\mathfrak{u}_{S}$ 
of $\mathfrak{q}_{S}$ is therefore
$\bigoplus_{\alpha\in R_{S}^{1}} \mathfrak{g}_{\alpha}$
where $R_{S}^{1}=R^{+}\setminus R_{S}$.

Let $S\subset \Pi$ be such that the unipotent radical $\mathfrak{u}_{S}$ 
of $\mathfrak{q}_{S}$ is abelian. Then it is easy to check that 
$\mathfrak{q}_{S}$ is maximal. Moreover, if we set $\mathfrak{u}_{S}^{-}=
\bigoplus_{\alpha\in R_{S}^{1}} \mathfrak{g}_{-\alpha}$, the decomposition
$$
\mathfrak{g}=\mathfrak{u}_{S}\oplus (\mathfrak{h}\oplus \bigoplus_{\alpha\in
R_{S}}\mathfrak{g}_{\alpha} )\oplus \mathfrak{u}_{S}^{-}
$$
induces a $\mathbb{Z}$-grading on $\mathfrak{g}$. This in turn 
defines a symmetric pair with 
$$
\mathfrak{k}_{S}=\mathfrak{h}\oplus \bigoplus_{\alpha\in
R_{S}}\mathfrak{g}_{\alpha}
\ , \
\mathfrak{p}_{S}=\mathfrak{u}_{S}\oplus \mathfrak{u}_{S}^{-}.
$$
These symmetric pairs correspond to inner automorphisms of
$\mathfrak{g}$ with non semisimple fixed point sets. 
A list of all such parabolic subalgebras and the corresponding 
symmetric pairs is given in Table 1 (see also Table 7 of the 
Reference chapter of~\cite{OV}).

We shall construct a nice Cartan subspace for these symmetric pairs
by using the following properties of the root system, sometimes
referred to as the ``cascade'' construction of Kostant. We shall
recall this construction and certain useful properties.
For more details, the reader may refer to \cite{JA}, \cite{JO} or \cite{TY1}.

Let $T\subset \Pi$. If $R_{T}$ is irreducible, then $T$ is connected,
and we shall denote by $\varepsilon_{T}$ the highest root of $R_{T}$.

We define a set $\mathcal{K}(T)$ by induction on the cardinality
of $T$ as follows :
\begin{itemize}
\item $\mathcal{K}(\emptyset )=\emptyset$.
\item If $T_{1},\dots ,T_{r}$ are the connected components of $T$,
then :
$$
\mathcal{K}(T)=\mathcal{K}(T_{1})\cup \cdots \cup \mathcal{K}(T_{r}).
$$
\item If $T$ is connected, then 
$$
\mathcal{K}(T)=\{ T\} \cup 
\mathcal{K}\bigl( \{ \alpha\in T \ ; \ \alpha \ \mathrm{and} \
\varepsilon_{T} \mbox{ are orthogonal}\} \bigr). 
$$ 
\end{itemize}
Note that any $K\in \mathcal{K}(T)$ is connected, and the roots
$\varepsilon_{K}$, $K\in \mathcal{K}(T)$, are pairwise 
strongly orthogonal (two distinct positive roots $\alpha,\beta$ are 
{\it strongly orthogonal}
if $\alpha\pm\beta$ are not roots).

For $K\in \mathcal{K}(T)$, denote by $\Gamma^{K}$
the set of positive roots $\alpha\in R_{K}^{+}$ which are not
orthogonal to $\varepsilon_{K}$. 

We have the following properties of the sets $\Gamma^{K}$.
\begin{itemize}
\item $R^{+}_{T}=\bigcup_{K\in \mathcal{K}(T )} \Gamma^{K}$ (disjoint union).
\item Let $K\in \mathcal{K}(T)$.
If $\alpha\in \Gamma^{K}\setminus \{ \varepsilon_{K}\}$, 
then there exists a unique $\beta\in \Gamma^{K}$ such that
$\alpha+\beta=\varepsilon_{K}$.
\item Let $K,K'\in \mathcal{K}(T )$,
$\alpha\in \Gamma^{K}$ and $\beta\in \Gamma^{K'}$
be such that $\alpha+\beta\in R$, then either
$K\subset K'$ and $\alpha+\beta \in \Gamma^{K'}$ or
$K'\subset K$ and $\alpha+\beta \in \Gamma^{K}$.
\end{itemize}

Let $S\subset \Pi$ be such that $\mathfrak{u}_{S}$ is abelian. We
set 
$$
\mathcal{E}=\{ K\in \mathcal{K}(\Pi ) \ ; \ \varepsilon_{K}\in R_{S}^{1}\},
$$ 
and for $K\in \mathcal{E}$, $X_{K}=X_{\varepsilon_{K}}+X_{-\varepsilon_{K}}$.
Note that the elements $X_{K}$ are semisimple. Let $(\mathfrak{g},\mathfrak{k}_{S})$
be the corresponding symmetric pair as defined above.

\begin{proposition}\label{sec}
\begin{enumerate}
\item We have $R_{S}^{1}\subset \bigcup_{K\in \mathcal{E}}\Gamma^{K}$.
\item Let $K\in \mathcal{E}$ and $\alpha\in R_{S}^{1}\cap \Gamma^{K}$.
Then $\varepsilon_{K}-\alpha \not\in R_{S}^{1}$.
\item The subspace $\mathfrak{a}$ spanned by the $X_{K}$, $K\in \mathcal{E}$,
is a Cartan subspace of $\mathfrak{p}_{S}$, and the rank of 
$(\mathfrak{g},\mathfrak{k}_{S})$ is $\sharp \mathcal{E}$.
\end{enumerate}
\end{proposition}
\begin{proof}
The first two parts are easy consequences of the properties of $\Gamma^{K}$
stated above and of the fact that $\mathfrak{u}_{S}$ is an abelian ideal of 
$\mathfrak{q}_{S}$.

Since the $\varepsilon_{K}$, $K\in \mathcal{E}$, are pairwise strongly 
orthogonal, the subspace $\mathfrak{a}$ of $\mathfrak{p}_{S}$ consists of pairwise
commuting semisimple elements, and its dimension is $\sharp \mathcal{E}$.
Using Table 1, we check that $\mathfrak{a}$ has the right dimension for
a Cartan subspace.
\end{proof}

We list below all the standard parabolic subalgebras $\mathfrak{q}_{S}$
with abelian unipotent radical, the unique simple root in $R_{S}^{1}$,
the corresponding symmetric pair, and the elements of the set
$\mathcal{E}$. We follow the numbering of simple roots in \cite{BO}.

$$
\begin{array}{c}
\vbox{\offinterlineskip
\def\cc#1{\hfill\kern .3em#1\kern .3em\hfill}
\def\tv{\vrule height 12pt depth 5pt}
\halign{\tv\cc{#}&\tv\cc{#}&\tv\cc{#}&\tv#\cr
\noalign{\hrule}
 & symmetric pair & elements of $\mathcal{E}$ &\cr
\noalign{\hrule}
$(A_{n},\alpha_{i})$ & $(\mathrm{sl}_{n+1},\mathrm{sl}_{n+1-i}\times 
\mathrm{sl}_{i}\times \mathbb{C})$&
$\{ \alpha_{j},\dots ,\alpha_{n+1-j} \}$,&\cr 
$n\geq 1$ & & $1\leq j\leq \min (i,n+1-i)$&\cr
\noalign{\hrule}
$(B_{n},\alpha_{1})$ & $(\mathrm{so}_{2n+1}, 
\mathrm{so}_{2n-1}\times \mathrm{so}_{2})$ & 
$\{ \alpha_{1}\} , \Pi $&\cr
$n\geq 2$ & & &\cr
\noalign{\hrule}
$(C_{n},\alpha_{n})$ & $(\mathrm{sp}_{2n}, \mathrm{gl}_{n})$ &
$\mathcal{K}(\Pi )$&\cr
$n\geq 2$ & & &\cr
\noalign{\hrule}
$(D_{n},\alpha_{1})$ & $(\mathrm{so}_{2n}, 
\mathrm{so}_{2n-2} \times \mathrm{so}_{2})$ &
$\{ \alpha_{1} \}, \Pi $&\cr
$n\geq 4$ & & &\cr
\noalign{\hrule}
$(D_{n},\alpha_{i})$,  $n\geq 4$ & $(\mathrm{so}_{2n} , 
\mathrm{gl}_{n})$ &
$\{ \alpha_{m},\dots ,\alpha_{n}\}, m\leq n-2$ odd&\cr
$i=n-1,n$ & & and $\{\alpha_{i}\}$ if $n\in 2\mathbb{Z}$
&\cr
\noalign{\hrule}
$(E_{6},\alpha_{1}),(E_{6},\alpha_{6})$& $(E_{6}, D_{5} \times \mathbb{C})$ &
$\Pi\setminus \{ \alpha_{2} \}$, $\Pi$&\cr
\noalign{\hrule}
$(E_{7},\alpha_{7})$ & $(E_{7}, E_{6}\times \mathbb{C})$ &
$\{ \alpha_{7}\}$, $\Pi\setminus \{\alpha_{1}\}$, $\Pi$&\cr
\noalign{\hrule}
}}\\
\mbox{Table 1}\cr
\end{array}
$$

\section{Symmetric pair associated to the centralizer of a semisimple element}
\label{centralizers}

In this section, $(\mathfrak{g},\mathfrak{k}_{S})$
is a symmetric pair associated to a standard parabolic subalgebra
$\mathfrak{q}_{S}$ with abelian unipotent radical.
Let $\mathfrak{a}$ be the Cartan subspace in $\mathfrak{p}_{S}$
as defined in part 3 of Proposition~\ref{sec}. 

Let $X\in \mathfrak{a}$. Then $\mathfrak{g}^{X}$ is a Levi factor
of a parabolic subalgebra of $\mathfrak{g}$.  Denote by 
$\mathfrak{l}=[\mathfrak{g}^{X},\mathfrak{g}^{X}]$ the semisimple
part of $\mathfrak{g}^{X}$, and set
$\mathfrak{l}_{+}=\mathfrak{l}\cap \mathfrak{k}_{S}^{X}$,
$\mathfrak{l}_{-}=\mathfrak{l}\cap\mathfrak{p}_{S}^{X}$
and $\mathfrak{r}_{+}=[\mathfrak{l}_{-},\mathfrak{l}_{-}]$. Then 
the decompositions
$$
\mathfrak{g}^{X}=\mathfrak{k}^{X}\oplus \mathfrak{p}^{X} \ , \
\mathfrak{l}=\mathfrak{l}_{+} \oplus \mathfrak{l}_{-} \ \mathrm{and} \
\mathfrak{r}=\mathfrak{r}_{+} \oplus \mathfrak{l}_{-}
$$
define symmetric subpairs of $(\mathfrak{g},\mathfrak{k}_{S})$, 
and the ranks of the pairs $(\mathfrak{l},\mathfrak{l}_{+})$ and  
$(\mathfrak{r},\mathfrak{r}_{+})$
are strictly inferior to that of $(\mathfrak{g},\mathfrak{k}_{S})$.

We shall determine, in certain cases, the symmetric pair
$(\mathfrak{r},\mathfrak{r}_{+})$ for any non-zero non $\mathfrak{p}_{S}$-regular
element $X\in\mathfrak{a}$, {\it i.e.} $\mathfrak{p}_{S}^{X}$ contains
a non-zero nilpotent element.

To simplify notations, $B_{0}$ will correspond to the Lie algebra $\{ 0\}$, 
$B_{1}=A_{1}$, $D_{2}=A_{1}\times A_{1}$ and $D_{3}=A_{3}$.

We shall need the classification of simple symmetric pairs of rank less than
or equal to $2$.

\begin{list1}\label{class1}
Rank $1$ :
$$
\begin{array}{ccc} 
(\mathrm{sl}_{n+1},\mathrm{sl}_{n}\times \mathbb{C}) , &
(\mathrm{so}_{n+1},\mathrm{so}_{n}) , \\
(\mathrm{sp}_{2n}, \mathrm{sp}_{2n-1}\times\mathrm{sp}_{2}) , & 
(F_{4},B_{4}) .
\end{array}
$$
\end{list1}
\begin{list2}\label{class2}
Rank $2$ :
$$
\begin{array}{cccc} 
(\mathrm{sl}_{n+2},\mathrm{sl}_{n}\times\mathrm{sl}_{2}\times \mathbb{C} ) , &
(\mathrm{sl}_{3} , \mathrm{so}_{3} ) , &
(\mathrm{sl}_{6} , \mathrm{sp}_{6} ) , &
(\mathrm{so}_{n+2}, \mathrm{so}_{n}\times \mathrm{so}_{2} ) , \\ 
(\mathrm{sp}_{2n+4},\mathrm{sp}_{2n}\times \mathrm{sp}_{4}) \ , &
(\mathrm{sp}_{4} , \mathrm{sl}_{2} \times \mathbb{C} ) \ , &
(\mathrm{so}_{10} , \mathrm{sl}_{5} \times \mathbb{C} ) \ , &
(E_{6} , F_{4} ) \ ,  \\ 
& (E_{6} , D_{5} \times \mathbb{C} ) \ , &
(G_{2} , A_{1}\times A_{1} ).
\end{array}
$$
\end{list2}

\begin{lemma}\label{typeB}
Suppose that we are in the case $(B_{n},\alpha_{1})$ of Table $1$. Then
there exists $m\in \mathbb{N}$ such that
$(\mathfrak{r},\mathfrak{r}_{+})=(\mathrm{so}_{m+1},\mathrm{so}_{m})$.
\end{lemma}
\begin{proof}
From the definition of $\mathfrak{a}$, we verify easily that $\mathfrak{a}$
commutes with the Lie subalgebra $\mathfrak{s}$ 
generated by the root vectors $X_{\pm\alpha}$, 
$\alpha \in \Pi \setminus \{ \alpha_{1},\alpha_{2} \}$.
So $\mathfrak{k}_{S}^{X}$ contains $\mathfrak{s}$.
Note that $\mathfrak{s}$ is simple of type $B_{n-2}$.

If $\mathfrak{t}$ is a Cartan subalgebra of the $\mathfrak{s}$,
then $\mathfrak{a}\oplus \mathfrak{t}$ is a Cartan subalgebra of 
$\mathfrak{g}^{X}$. It follows that the root system of the 
semisimple part $\mathfrak{l}$ of $\mathfrak{g}^{X}$ contains as
a subsystem the root system of $\mathfrak{s}$. In particular,
the semisimple rank of $\mathfrak{l}$ is equal to $n-2$ or $n-1$.

Since $X$ is not $\mathfrak{p}$-regular, $\mathfrak{p}^{X}$ 
contains a nilpotent element, and so $\mathfrak{l}$ contains
strictly $\mathfrak{s}$. 
It follows that the root system $R(\mathfrak{l})$ of $\mathfrak{l}$ 
is of one of the following types : 
$A_{n-1}$, $B_{n-1}$ or $A_{1}\times B_{n-2}$. 

Suppose that $R(\mathfrak{l})$ is not of type $A_{n-1}$. 
Then $\mathfrak{l}_{+}$ contains a simple Lie subalgebra of type
$B_{n-2}$. Consequently, using the classification of rank $1$ 
symmetric pairs (see List \ref{class1}), 
we deduce that $\mathfrak{r}$ has the required form.

Suppose now that $R(\mathfrak{l})$ is of type $A_{n-1}$. Then 
from the fact that $R(\mathfrak{l})$ contains 
a subsystem of type $B_{n-2}$, we deduce that $n=3$.
Explicit computations as described in the Section~\ref{computations} 
show that this case does not occur.
\end{proof}

\begin{lemma}\label{typeD}
Suppose that we are in the case $(D_{n},\alpha_{1})$ of Table $1$. Then
there exists $m\in \mathbb{N}$ such that
$(\mathfrak{r},\mathfrak{r}_{+})=(\mathrm{so}_{m+1},\mathrm{so}_{m})$.
\end{lemma}
\begin{proof}
The proof is more or less the same as in the case $(B_{n},\alpha_{1})$.
However, note that $\mathfrak{s}$ is simple except when $n=4$,
in which case it is semisimple. The root system $R(\mathfrak{l})$
is of one of the following types : 
$A_{n-1}$, $D_{n-1}$ or $A_{1}\times D_{n-2}$.
Here, if $R(\mathfrak{l})$ is of type $A_{n-1}$, then
$n=4$ or $n=5$. In the former, we use the classification of rank $1$
symmetric pairs to obtain a contradication, while in the latter,
a direct computation is required (see Section~\ref{computations}).
\end{proof}

\begin{remark}\label{smallrank}\rm
Note that we may extend Lemma~\ref{typeD}
to the symmetric pairs 
$(\mathrm{so}_{4},\mathrm{so}_{2}\times \mathrm{so}_{2})$ and
$(\mathrm{so}_{6},\mathrm{so}_{4}\times \mathrm{so}_{2})$.

In the first case, we have 
$(\mathrm{so}_{4},\mathrm{so}_{2}\times \mathrm{so}_{2})
=(\mathrm{so}_{3},\mathrm{so}_{2})\times
(\mathrm{so}_{3},\mathrm{so}_{2})$, while in the second case,
we have
$(\mathrm{so}_{6},\mathrm{so}_{4}\times \mathrm{so}_{2})
=(\mathrm{sl}_{4},\mathrm{sl}_{2}\times \mathrm{sl}_{2}\times \mathbb{C} )$.
In both cases, we may obtain the same conclusion by direct computations
(see Section~\ref{computations}).
\end{remark}

Summarizing, since Cartan subspaces are $K$-conjugate,
we have therefore obtained the following result :

\begin{proposition}\label{so}
Let $(\mathfrak{g},\mathfrak{k})$ be the symmetric pair 
$(\mathrm{so}_{p+2},\mathrm{so}_{p}\times \mathrm{so}_{2})$,
$p\geq 2$. For any non-zero non $\mathfrak{p}$-regular semisimple
element $X$ in $\mathfrak{p}$, the subpair $(\mathfrak{r},\mathfrak{r}_{+})$
is of the form $(\mathrm{so}_{m+1},\mathrm{so}_{m})$ for some $m\in\mathbb{N}$.
\end{proposition}

\begin{lemma}\label{typeE7}
Suppose that we are in the case $(E_{7},\alpha_{7})$ of Table 1. Then
the symmetric pair $(\mathfrak{r},\mathfrak{r}_{+})$ is the product
of symmetric pairs of the form
$(\mathrm{so}_{m+1},\mathrm{so}_{m})$
or $(\mathrm{so}_{m+2},\mathrm{so}_{m}\times \mathrm{so}_{2})$ 
for some integer $m\in\mathbb{N}$.
\end{lemma}
\begin{proof}
Proceeding as in the proof of Lemma~\ref{typeB}, we check easily
that $\mathfrak{k}_{S}^{X}$ contains the Lie subalgebra $\mathfrak{s}$ 
generated by the root vectors $X_{\pm\alpha}$, 
$\alpha \in \Pi \setminus \{ \alpha_{1},\alpha_{6},\alpha_{7} \}$.
Note that $\mathfrak{s}$ is simple of type $D_{4}$.

Hence $R(\mathfrak{l})$ contains a subsystem of type $D_{4}$. 
We deduce therefore that $R(\mathfrak{l})$ is of one of the following
types : $D_{5}$, $D_{4}\times \mathbb{C}$, $D_{6}$ or $D_{5}\times A_{1}$.
We may conclude by using the classification of symmetric pairs of
rank less than or equal to $2$ (see Lists \ref{class1} and \ref{class2}).
\end{proof}

\begin{remark}\label{rank2}\rm
Note that in all the other rank $2$ cases, there exists 
$X\in \mathfrak{a}$ non $\mathfrak{p}$-regular such that 
$\mathfrak{p}^{X}$ contains
two non proportional commuting nilpotent elements. Hence
by the result of \cite[Proposition 3]{SY}, the conclusion of the previous lemmas
fails. See Section~\ref{computations} for more details.
\end{remark}

\section{Sheets and commuting varieties}\label{commuting}

Let $(\mathfrak{g},\mathfrak{k})$ be a symmetric pair.
Recall that the connected algebraic group $K$ acts on $\mathfrak{p}$.
For $n\in\mathbb{N}$, we set :
$$
\mathfrak{p}^{(n)}=\{ X\in \mathfrak{p} \ ; \ \dim K.X=n \}.
$$
The set $\mathfrak{p}^{(n)}$ is locally closed, and an irreducible component
of $\mathfrak{p}^{(n)}$ shall be called a {\em $K$-sheet} of $\mathfrak{p}$.
Clearly, $K$-sheets are $K$-invariant, and 
by \cite{T}, each $K$-sheet contains a nilpotent element.

Let $\pi_{1} : C(\mathfrak{p})\rightarrow \mathfrak{p}$
be the projection $(X,Y)\mapsto X$.
Recall the following result concerning the commuting
variety of $\mathfrak{p}$.

\begin{theorem}\label{sheet}
There exist $K$-sheets $\mathcal{S}_{1},\dots , \mathcal{S}_{r}$ of
$\mathfrak{p}$
such that $\overline{\pi_{1}^{-1}(\mathcal{S}_{i})}$, $i=1,\dots ,r$,
are the irreducible components of $C(\mathfrak{p})$.
\end{theorem}

The proof of Theorem~\ref{sheet} is a simple consequence of the
following result. For the sake of completeness, 
we have included a proof.

\begin{lemma}\label{fiber}
Let $V$ be a vector space, $E\subset V\times V$ a locally closed subvariety
and for $i=1,2$, $\pi_{i}:E\rightarrow V$ be 
the projection $(x_{1},x_{2})\mapsto x_{i}$. 
Suppose that :
\begin{itemize}
\item[1.] $\pi_{1}(E)$ is locally closed.
\item[2.] There exists $r\in\mathbb{N}$ such that
for all $x\in \pi_{1}(E)$,  $\pi_{2}(\pi_{1}^{-1}(x))$ is a vector 
subspace of dimension $r$. 
\end{itemize}
If $\pi_{1}(E)$ is irreducible, then so is $E$.
\end{lemma} 
\begin{proof}
Let $\mathbf{G}$ be the Grassmann variety of 
$r$-dimensional subspaces of $V$, $x\in \pi_{1}(E)$ and 
$W=\pi_{2} (\pi_{1}^{-1}(x))\in \mathbf{G}$.
Fix a complementary subspace $U$ of $W$ in $V$ and set :
$$
\mathbf{F} =\{ T \in \mathbf{G} \ ; \ T \cap U = \{ 0\} \}
=\{ T \in \mathbf{G} \ ; \ T + U = V \}.
$$
Clearly, $\mathbf{F}$ is an open subset of $\mathbf{G}$ containing 
$W$. For $\tau \in \Hom (W,U)$ the set of linear maps from $W$
to $U$, we define:
$$
T(\tau ) = \{ w +\tau (w) \ ; \ w\in W \}.
$$
Then we check easily that $T( \tau )\in \mathbf{F}$,
and we have a map
$
\Hom (W,U) \rightarrow \mathbf{F} \ , \ \tau \mapsto T( \tau ).
$ 
We claim that this map is an isomorphism.

Since $w_{1}+\tau_{1}(w_{1})=w_{2}+\tau_{2}(w_{2})$
is equivalent to
$w_{1}-w_{2}=\tau_{2}(w_{2})-\tau_{1}(w_{1})$,
we deduce that the above map is injective.

Now if $T\in \mathbf{F}$, then 
for $w\in W$, we define $\tau (w)$ to be the unique 
element in $U$ such that $w+\tau (w)\in T$.
We then verify easily that $T(\tau )= T$.
So we have proved our claim.

The map
$$
\Phi : \pi_{1}(E ) \rightarrow \mathbf{G} \ , \ 
y\mapsto \pi_{2} (\pi_{1}^{-1}(y)) 
$$
is a morphism of algebraic varieties.
So $F=\Phi^{-1}(\mathbf{F})$ is an open subset of $\pi_{1}(E)$
containing $x$. The above claim says that we have a well-defined map:
$$
\Psi : F \times W \rightarrow E \ , \
(y,w) \mapsto (y,w+\tau (w) )
$$
where $T( \tau )=\Phi (y)$. It is then a straightforward verification
that $\Psi$ is an isomorphism of the algebraic varieties
$F\times W$ and $\pi_{1}^{-1}(F)$.

It follows that the map
$\pi_{1} : E \rightarrow \pi_{1}(E)$ 
is an open map whose fibers are irreducible. 
Hence by a classical result on topology \cite[T.5]{DI}, 
if $\pi_{1}(E)$ is irreducible, then $E$ is irreducible.
\end{proof}

Since the set of $\mathfrak{p}$-generic elements and the
set $\mathfrak{p}_{\mathrm{reg}}$ of $\mathfrak{p}$-regular 
elements are open subsets of $\mathfrak{p}$, we have the following corollary :

\begin{corollary}\label{maximal}
Let $\mathfrak{a}$ be a Cartan subspace in $\mathfrak{p}$.
The set 
$$
C_{0}=\overline{K.(\mathfrak{a}\times \mathfrak{a})}
=\overline{\pi_{1}^{-1}(\mathfrak{p}_{\mathrm{reg}})}
=\overline{\pi_{2}^{-1}(\mathfrak{p}_{\mathrm{reg}})}
$$
is the unique irreducible component of $C(\mathfrak{p})$ of
maximal dimension. 
\end{corollary}

We shall finish this section with a result on even nilpotent
elements. Let $X\in \mathfrak{p}$ be a nilpotent element,
and $(H,Y)\in \mathfrak{k}\times \mathfrak{p}$ be such that
$(X,H,Y)$ is a normal $\mathrm{sl}_{2}$-triple
(called a normal $\mathrm{S}$-triple in \cite{T}).
Recall that $X$ is {\em even} if the eigenvalues of 
$\ad_{\mathfrak{g}} H$ are even.  In fact, this is equivalent
to the condition that the eigenvalues of 
$\ad_{\mathfrak{p}} H$ are even.

\begin{proposition}\label{even}
Let $X\in \mathfrak{p}$ be an even nilpotent element, then  
$X$ belongs to a $K$-sheet containing semisimple elements.
\end{proposition}
\begin{proof}
Let $(X,H,Y)$ be a normal $\mathrm{sl}_{2}$-triple and
$\mathfrak{s}=\mathbb{C} X+\mathbb{C} H+\mathbb{C} Y$. 
Then $\mathfrak{g}$ decomposes into a direct sum
of $\mathfrak{s}$-modules, say $V_{i}$, $i=1,\dots ,r$. Since $X$ is even, 
$\dim V_{i}$ is odd for $i=1,\dots ,r$. 

For $\lambda \in \mathbb{C}$, we set $X_{\lambda} = X+\lambda Y\in \mathfrak{p}$.
If $\lambda\neq 0$, then $X_{\lambda}$ is semisimple
because $X_{\lambda}$ is $G$-conjugate to a multiple of $H$. 
We claim that $\dim \mathfrak{p}^{X_{\lambda}}=\dim \mathfrak{p}^{X}$
for all $\lambda \in \mathbb{C}$.

First of all, observe that $\mathfrak{p}^{X_{\lambda}}
=\bigoplus_{i=1}^{r} (V_{i}\cap \mathfrak{p})^{X_{\lambda}}$ because
$V_{i}=(V_{i}\cap \mathfrak{k})\oplus (V_{i}\cap \mathfrak{p})$.
Moreover $\dim (V_{i}\cap \mathfrak{p})^{X_{\lambda}}\leq 1$.

Now if $(V_{i}\cap \mathfrak{p})^{X_{\lambda}}\neq \{ 0\}$, then 
a simple weight argument shows that $(V_{i}\cap \mathfrak{p})^{X}\neq \{ 0\}$.

Conversely, suppose that $(V_{i}\cap \mathfrak{p})^{X}\neq \{ 0\}$.   
Let $\dim V_{i}=2n+1$ and 
$v_{-n}, \dots ,v_{n}$ be a basis weight vectors of $V_{i}$ such that 
$H v_{k}=2k v_{k}$, $k=-n,\dots ,n$. 
Then $(V_{i}\cap \mathfrak{p})^{X}=\mathbb{C} v_{n}$. 
 
So $v_{k}\in \mathfrak{k}$ (resp. $v_{k} \in \mathfrak{p}$) when $n-k$ is odd
(resp. even). In particular, $v_{-n} \in \mathfrak{p}$.
It follows that for $k$ such that $n-k$ odd,
$\lambda Y v_{k+1} =  - a_{k} X v_{k-1}$ for some $a_{k}\in \mathbb{C}$.
We may therefore renormalize the $v_{k}$'s so that
$
v = v_{-n} + v_{-n+2}+\cdots + v_{n-2} +v_{n}
$
verifies $X_{\lambda} v=0$.

We have therefore proved that $\dim \mathfrak{p}^{X_{\lambda}}=\dim \mathfrak{p}^{X}$
for all $\lambda$.

Now, consider the morphism 
$\Phi : K \times \mathbb{C}  \rightarrow \mathfrak{p}$, $(k,\lambda ) \mapsto
k.X_{\lambda}$. The image of $\Phi$ is irreducible and 
contains semisimple elements, so it contains strictly
$K.X$. Consequently, $K.X$ is contained strictly in a $K$-sheet with semisimple
elements.
\end{proof}

\section{The rank 2 case}\label{ranktwo}

In this section, $(\mathfrak{g},\mathfrak{k})$ is a symmetric pair of rank
$2$.

\begin{lemma}\label{semisimple}
Let $X \in \mathfrak{p}$ be an element which is non-semisimple and non-nilpotent, 
and $X_{s}$, $X_{n} \in \mathfrak{p}\setminus \{0\}$ be its semisimple
and nilpotent components. If $(X,Y) \in C(\mathfrak{p})$, then the
semisimple component of $Y$ belongs to $\mathbb{C} X_{s}$. 
\end{lemma}
\begin{proof}
Let $Y = Y_{s} + Y_{n}$ be the corresponding decomposition into semisimple
and nilpotent components. Since
$\mathfrak{g}^X = \mathfrak{g}^{X_{s}} \cap \mathfrak{g}^{X_{n}}$, 
it follows easily that $\mathfrak{p}^{X} = \mathfrak{p}^{X_{s}} \cap 
\mathfrak{p}^{X_{n}}$. Hence
$X_{s}$, $X_{n}$, $Y_{s}$, $Y_{n}$ commute pairwise. If $X_{s}$ and $Y_{s}$ are not
proportional, then they span a Cartan subspace because the rank is
two. But this would imply that $X$ commutes with a Cartan subspace,
thus $X$ belongs to a Cartan subspace, which is absurd because $X$ is
not semisimple.
\end{proof}

In the rest of this section, we consider the
rank $2$ symmetric pair 
$(\mathfrak{g},\mathfrak{k})= (\mathrm{so}_{p+2}, \mathrm{so}_{p}\times
\mathrm{so}_{2})$, $p\geq 2$. 
In particular, if $p \neq 2$, then the symmetric pair 
comes from a parabolic subalgebra with abelian unipotent radical.

Let $\mathfrak{a}$ be a Cartan subspace in $\mathfrak{p}$. 
Recall that $C_{0} = \overline{K.(\mathfrak{a} \times \mathfrak{a})}$.

\begin{theorem}\label{reduction}
The commuting variety $C(\mathfrak{p})$ is irreducible.
\end{theorem}
\begin{proof}
We proceed as in the proof of Richardson in the classical case (see \cite{R}) 
by using inductive arguments. Let $(X,Y) \in C(\mathfrak{p})$.
\vskip0.5em

1. If $X$ is semisimple, then $X$ commutes with a $\mathfrak{p}$-regular
semisimple element $Z$. The line $\mathcal{L}_{Z} = \{ (X,tY + (1-t)Z), t \in \mathbb{C} \}$ 
is contained in $C(\mathfrak{p})$. Since $\{ tY + (1-t)Z, t \in \mathbb{C}\}$ meets the set of 
$\mathfrak{p}$-regular semisimple elements which is open in $\mathfrak{p}$, we
conclude that ${\mathcal L}_{Z}$, and hence $(X,Y)$, are contained in $C_{0}$
(Corollary~\ref{maximal}).
\vskip0.5em

2. We may assume that neither $X$ nor $Y$ is semisimple.

Suppose that $X$ is not nilpotent. Writing $X = X_{s} + X_{n}$, $Y = Y_{s} +
Y_{n}$, we deduce from the Lemma~\ref{semisimple} that $Y_{s} \in \mathbb{C} X_{s},
X_{n},Y_{n} \in \mathfrak{p}^{X_{s}}$ and $[X_{n},Y_{n}] = 0$.  

By Proposition~\ref{so}, we deduce that
$X_{n},Y_{n}$ are commuting elements contained in a symmetric pair of the form
$(\mathrm{so}_{m+1},\mathrm{so}_{m})$ for some $m\in \mathbb{N}$. Consequently, 
$Y_{n} \in \mathbb{C} X_{n}$ by \cite[Proposition 3]{SY}.

Thus, there exists $\lambda,\mu \in \mathbb{C}$ such that $Y =
\lambda X_{s} + \mu X_{n}$. In particular, $X$ is $\mathfrak{p}$-regular, and
therefore $(X,Y) \in C_{0}$ by Corollary~\ref{maximal}.
\vskip0.5em

3. So we may further assume that $X$ and $Y$ are both nilpotent. If $X$
commutes with a non-zero semisimple element $Z \in \mathfrak{p}$, then
the same argument as in 1) works because the set of non-nilpotent
elements is open.
\vskip0.5em

4. Recall that an element of $\mathfrak{p}$ is said to be 
$\mathfrak{p}$-distinguished if its centralizer in $\mathfrak{p}$ 
does not contain any non-zero semisimple element. In particular, 
a $\mathfrak{p}$-distinguished element is nilpotent. So the number
of $K$-orbits of $\mathfrak{p}$-distinguished elements is finite.

Denote by $\pi_{1} : C(\mathfrak{p})\rightarrow \mathfrak{p}$ 
the projection $(X_{1},X_{2})\mapsto X_{1}$,
$\mathcal{O}$ the set of non $\mathfrak{p}$-distinguished elements
in $\mathfrak{p}$, and $\Omega_{1},\dots ,\Omega_{r}$ the set
of $K$-orbits of $\mathfrak{p}$-distinguished elements in $\mathfrak{p}$.
Thus 
$\mathfrak{p}=\mathcal{O} \cup \Omega_{1}\cup \cdots \cup \Omega_{r}$,
and
$C(\mathfrak{p})=\pi_{1}^{-1}(\mathcal{O})\cup
\pi_{1}^{-1}(\Omega_{1})\cup \cdots \cup \pi_{1}^{-1}(\Omega_{r})$.

From the previous paragraph, we obtain that
$\pi^{-1}_{1}(\mathcal{O})\subset C_{0}$.
Consequently, $C(\mathfrak{p})$ is the union of $C_{0}$ with
$\overline{\pi_{1}^{-1}(\Omega_{i_{1}})}, \dots, 
\overline{\pi_{1}^{-1}(\Omega_{i_{s}})}$. 
Now we check easily that for $X\in\mathfrak{p}$,
$\pi_{1}^{-1}(K.X) = K.(X,\mathfrak{p}^{X})$ is an irreducible subset of
$C(\mathfrak{p})$ of dimension $\dim \mathfrak{k} - \dim \mathfrak{k}^{X} + 
\dim \mathfrak{p}^{X} = \dim \mathfrak{p}$. 

It follows that all irreducible components of $C(\mathfrak{p})$ other than
$C_{0}$, if they exist, are of dimension $\dim \mathfrak{p}$.

Suppose that $C(\mathfrak{p})$ is not irreducible.
By the previous discussion, there exists a $\mathfrak{p}$-distinguished
element $X$ such that
$\overline{\pi_{1}^{-1}(K.X)}$ is an irreducible component of 
dimension $\dim \mathfrak{p}$.

If $X$ is even, then by Proposition~\ref{even},
$X$ belongs to a $K$-sheet $\mathcal{S}$ 
containing non-zero semisimple elements. So $\dim \mathcal{S} > \dim K.X$.
Now Lemma~\ref{fiber} says that
$\dim \pi_{1}^{-1} (\mathcal{S}) > \dim \mathfrak{p}$. Contradiction.

So we may assume that $X$ is not even and the theorem follows
from the following result.
\end{proof}

\begin{proposition}\label{evenness}
All $\mathfrak{p}$-distinguished elements of 
the symmetric $(\mathrm{so}_{p+2},\mathrm{so}_{p}\times \mathrm{so}_{2})$,
$p\geq 2$, are even.
\end{proposition}

We shall first classify the nilpotent
$K$-orbits in $\mathfrak{p}$.

$\bullet$
Each element $X$ of $\mathfrak{g}$ can be represented by a matrix of the
following type :
$$
X = \left(
\begin{array}{cc}
A_{p} & B \\
-{}^tB & A_2 
\end{array}
\right), \ 
A_p \in \mathrm{so}_{p} , A_2 \in \mathrm{so}_{2}, B
\in \mathrm{M}_{p,2}(\mathbb{C} ).
$$

The involutive automorphism $\theta$ is defined by :

$$
\theta(X) = J_{p,2}XJ_{p,2}, \mbox{where } 
J_{p,2} = \left(
\begin{array}{cc}
I_p & 0 \\
0& -I_2 
\end{array}\right).
$$
It follows that :  
$$
\mathfrak{p} = \left\{ 
\left(\begin{array}{cc}
0 & B \\ 
-{}^tB& 0
\end{array}\right), \ B \in \mathrm{M}_{p,2}(\mathbb{C} ) \right\}.
$$
The dimension of $\mathfrak{p}$ is therefore $2p$.

Set :
$$
H_{i} = E_{i, n-i+1} - E_{n-i+1,i}\ ,\ 1 \leq i \leq 2.
$$ 
The subspace $\mathfrak{a}$, spanned by $H_{1},H_{2}$, is a Cartan
subspace of $\mathfrak{p}$.

$\bullet $
Let $\mathcal{N}$ be the set of $G$-nilpotent orbits in
$\mathfrak{g}$, $\mathrm{YD}_{p+2}$ the set of all Young diagrams 
corresponding to the partitions of $p+2$ and satisfying the two following 
properties :

\begin{itemize} 
\item[$\bf (P_1)$] Each row of even length occurs with even multiplicity.
\item[$\bf (P_2)$] If all rows have even length, then Roman numerals are
attached to the diagram, namely I or II.
\end{itemize}

Recall (see for example~\cite{CM}) that $\mathcal{N}$ is in one-to-one 
correspondence with $\mathrm{YD}_{p+2}$.

$\bullet $ Let $\mathcal{N}_{\mathfrak{p}}$ be the set of nilpotent
$K$-orbits in $\mathfrak{p}$.  To describe this set, we use the Kostant-Sekiguchi
correspondence between $\mathcal{N}_{\mathfrak{p}}$ and the set
$\mathcal{N}_{0}$ of nilpotent $G_{0}$-orbits in $\mathfrak{g}_{0}$,
where $\mathfrak{g}_{0}$ is a real form of $\mathfrak{g}$ and $G_{0}$ the
corresponding connected real Lie group. Let us recall briefly this correspondence. 

We need to fix some notations : 
\begin{itemize}
\item[-] $\mathfrak{g}_{0} = \mathrm{so}(p,2)$  is the set of  matrices $X_{0}$ of the
following type :
$$
X_{0} = \left(
\begin{array}{cc}
A_p & B \\ 
{}^tB& A_2 
\end{array}\right), \ 
A_p \in \mathrm{so}_{p}(\mathbb{R}), \ A_2 \in \mathrm{so}_{2}(\mathbb{R}), \ B
\in \mathrm{M}_{p,2}(\mathbb{R}).
$$
\item[-] A Cartan decomposition of $\mathfrak{g}_{0}$ is given by :
$$
\mathfrak{g}_{0} = \mathfrak{k}_{0} \oplus \mathfrak{p}_{0}
$$
where the Cartan involution $\theta_{0}$ is defined by :
$\theta_{0}(X_{0}) = - \ {}^tX_{0}$. 
$$
\begin{array}{rl}
\mathfrak{k}_{0} &= \left\{ 
\left(\begin{array}{cc}
A_p & 0 \\ 0& A_2 
\end{array}\right), \ 
A_p \in \mathrm{so}_{p}(\mathbb{R}), A_2 \in \mathrm{so}_{2}(\mathbb{R})\right\} ,\\
[0.2em]
\mathfrak{p}_{0} &= \left\{
\left(\begin{array}{cc}
0 & B \\ {}^tB& 0 
\end{array}\right), 
B \in \mathrm{M}_{p,2}(\mathbb{R}) \right\}.
\end{array}
$$
\item[-] $G_{0}$ is the adjoint group of $\mathfrak{g}_0$.
\end{itemize}

The Lie algebra $\mathfrak{g}_{0}$ is embedded in $\mathfrak{g}$, via the map 
$\phi$ defined by:
$$
\phi(X_0) = \widetilde{J}_{p,2}X_0\widetilde{J}_{p,2}^{-1}
$$
where :
$$
\widetilde{J}_{p,2} = \left(\begin{array}{cc}
I_p & 0 \\ 0& -iI_2
\end{array}\right)
$$ 

We check easily that : $\phi \circ \theta_{0} = \theta
\circ \phi$. It follows that $\mathfrak{g}_0$ can be identified with a real
form of $\mathfrak{g}$. Hence, we have :
$$
(\mathfrak{k}_{0})_{\mathbb{C}} = \mathfrak{k}, \  
(\mathfrak{p}_{0})_{\mathbb{C}} = \mathfrak{p}.
$$

$\bullet $ The next step is to describe ${\mathcal N}_{0}$, 
using the method described in~\cite[Chapter 9]{CM}.

First of all, we define a signed Young diagram of signature $(p,2)$  
to be a Young diagram in which each box is labeled with a $+$ or a $-$ sign 
in such a way that :
\begin{itemize}
\item[-] Signs alternate across rows.
\item[-] Two diagrams are equivalent if and only if one can be obtained from
the other by interchanging rows of equal length.
\item[-] The number of boxes labeled $+$ is $p$, the number of boxes labeled $-$
is $2$.
\end{itemize}

We define an orthogonal signed Young diagram of signature
$(p,2)$ to be a signed Young diagram satisfying $\bf (P_1)$, $\bf (P_2)$ and 
the three following properties :
\begin{itemize}
\item[$\bf (P_3)$] Each row of even length has its leftmost box labeled $+$.
\item[$\bf (P_4)$] If all rows have even length, then two Roman numerals,
each I or II, are attached to the corresponding diagram (giving
then four orbits).
\item[$\bf (P_5)$] If at least one row has odd length and if each row of
  odd length has an even number of boxes labeled $+$ or if each row of
  odd length has an even number of boxes labeled $-$, then one numeral I
  or II is attached to the corresponding diagram.
\end{itemize}

Let us denote by $\mathrm{DYO}_{p,2}$ the set of all orthogonal signed Young diagram.

By a classical result (see~\cite[9.3.4]{CM}), we know that $\mathcal{N}_{0}$ 
is in one-to-one correspondence with $\mathrm{DYO}_{p,2}$.

\begin{remark}\rm\label{rem1}
If $\Omega_{0} \in \mathcal{N}_0$, then $\Omega = 
(\Omega_{0})_{\mathbb{C} } \in \mathcal{N}$. 
The Young diagram corresponding to $\Omega$ is
obtained by removing the signs in the orthogonal signed Young diagram
corresponding to $\Omega_{0}$.
\end{remark}

\begin{remark}\rm\label{rem2}
It follows from the definition of elements of
$\mathrm{DYO}_{p,2}$ that if 
$\Omega_{0} \in \mathcal{N}_{0}$ and $\Omega_{0}= (d_1,\dots,d_s)$, then $d_{1} \leq 5$.
\end{remark}
 
$\bullet $ We may deduce now, from the previous statements, all the
nilpotent orbits in $\mathrm{so}(p,2)$.

If $p\geq 4$, we have :

$$
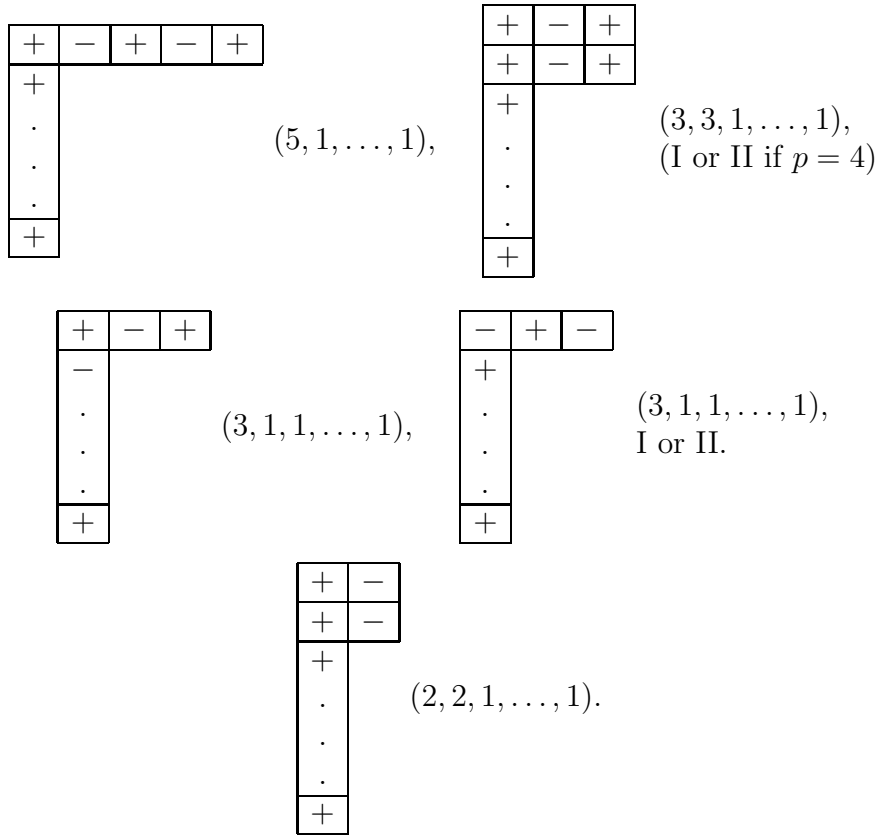
\begin{array}{|c|c|c|c|c|}
\hline 
+&-&+&-&+\\
\hline 
+\\
.\\
.\\ 
.\\
\cline{1-1} 
+\\
\cline{1-1}
\end{array}
\ 
(5,1,\dots,1),
\ \ \ \ 
\begin{array}{|c|c|c|}
\hline 
+ & - & +  \\
\hline 
+ & - & +  \\
\hline 
+ \\
.\\
.\\
.\\
\cline{1-1} +\\
\cline{1-1}
\end{array} 
\  
\begin{array}{l}
(3,3,1,\dots,1),\\ 
(\mbox{I or II if } p=4)
\end{array}
$$    
$$    
\begin{array}{|c|c|c|}
\hline + & - & +  \\
\hline - \\
.\\
.\\
.\\
\cline{1-1}+\\
\cline{1-1}
\end{array} 
\ 
(3,1,1,\dots,1),
\ \ \ \ 
\begin{array}{| c|c |c |}
\hline - & + & -  \\
\hline + \\
.\\
.\\
.\\
\cline{1-1} + \\
\cline{1-1}
\end{array}
\  
\begin{array}{l}
(3,1,1,\dots,1),\\
\mbox{I or II.} 
\end{array}
$$
$$
\begin{array}{|c|c|}
\hline + & -  \\
\hline + & - \\
\hline + \\
.\\
.\\
.\\
\cline{1-1} + \\
\cline{1-1}
\end{array} 
\ 
(2,2,1,\dots,1).
$$

If $p=3$, the orbits are given by the partitions $(5),(3,1), (2,2,1)$
and the corresponding diagrams  above.

If $p=2$, we have only the partitions $(3,1)$ and $(2,2)$.

$\bullet $ Let $\Omega_{0} \in \mathcal{N}_{0}$. We may choose a $\mathrm{sl}_{2}$-triple 
$(H_0,X_0,Y_0)$ such that :
 
$$
X_0 \in \Omega_{0}\ , \ \theta_0(H_0) = -H_0\ ,\  \theta_0(X_0) =
-Y_0\ ,\ \theta_0(Y_0) = -X_0.
$$
Triples with these properties are called {\it Cayley triples}.

Set :
$$
H_S = i(X_0 - Y_0)\ ,\ \ X_S = \frac{1}{2}(X_0 + Y_0 + iH_0)\ ,\ \ Y_S =
\frac{1}{2}(X_0 + Y_0 - iH_0)
$$
Then $(H_S,X_S,Y_S)$ is a normal $\mathrm{sl}_{2}$-triple (that is, $H_S
\in \mathfrak{k}$ and $X_S$, $Y_S \in \mathfrak{p}$), called the 
{\it Cayley transform} of $(H_0,X_0,Y_0)$.

\begin{theorem}[Kostant-Sekiguchi correspondence \cite{CM}]  
There is a bijection from  $\mathcal{N}_0$ onto  $\mathcal{N}_{\mathfrak{p}}$ 
which sends the orbit of the positive element of a Cayley triple to the orbit of the
positive element of its Cayley transform.
\end{theorem}

Via this correspondence, we classify the nilpotent $K$-orbits
in  $\mathfrak{p}$ using the set $\mathrm{DYO}_{p,2}$.

$\bullet $ Let  $r = \displaystyle \left[\frac{p+2}{2}\right]$ be the rank of 
$\mathfrak{g}$ and $\mathfrak{h} = \mathfrak{a} \oplus \mathfrak{t}$ a Cartan subalgebra of
$\mathfrak{g}$, where $\mathfrak{t}$ is a Cartan subalgebra of the centralizer of 
$\mathfrak{a}$ in $\mathfrak{k}$. We may complete
the basis $(H_{1},H_{2})$ of $\mathfrak{a}$ to a basis $(H_{i}), 1 \leq i \leq r$, 
of $\mathfrak{h}$, where $(H_i, 3 \leq i
\leq r)$,  is a basis of $\mathfrak{t}$. Let $(e_{i})_{1\leq i\leq r}$ 
be the corresponding 
dual basis in $\mathfrak{h}^{*}$, $R = R(\mathfrak{g}, \mathfrak{h})$ the 
corresponding root system and $\Pi$ a base of simple roots of $R$,  
defined as follows : 
$$
\begin{array}{rl} 
p+2 = 2r+1, \ & R = \{ \pm e_i \pm e_j, 1 \leq i \not= j \leq r\} \cup
\{\pm e_i, 1 \leq i \leq r\}\\
& \Pi = \{ e_i - e_j, 1 \leq i < j \leq
r\} \cup \{e_r\} \\   
p+2 = 2r,\ & R = \{ \pm e_i \pm e_j, 1 \leq i \not= j \leq r\} \\
&\Pi = \{ e_i - e_j , 1 \leq i < j \leq r\} \cup \{e_{r-1} + e_r\}
\end{array}
$$

In the same way, we choose a Cartan subalgebra $\mathfrak{h}_{0}$
of $\mathfrak{g}_{0}$, such that $\mathfrak{h} = (\mathfrak{h}_{0})_{\mathbb{C} }$. In
particular, $\mathfrak{a}_{0} = \mathbb{R} H_1 \oplus \mathbb{R}H_2$ is a
Cartan subspace of $\mathfrak{p}_{0}$ such that $\mathfrak{a} = (\mathfrak{a}_{0})_{\mathbb{C} }$.

$\bullet $ Given a nilpotent orbit $\Omega$ in $\mathcal{N}$,
there is a unique $\mathrm{sl}_{2}$-triple $(h^{+},e,f)$ such that :
$e \in \Omega$ and for all $\alpha_{i} \in \Pi$, $\alpha_i(h^{+}) \in
\{0,1,2\}$. 

Recall that the sequence $C(\Omega) = (\alpha_{1}(h^{+}),\dots ,\alpha_{r}(h^{+}) )$
is called the {\it characteristic} of $\Omega$.

Given a Young diagram, we can compute now the characteristic of the
corresponding nilpotent orbit $\Omega$ by using the following method :
\begin{itemize}
\item[-] First of all, we suppose that $p+2 = 2r+1$. Let  $\Omega = (d_1,
\dots,d_s)$ be a nilpotent $G$-orbit. We consider the sequence of integers 
$(d_i -1, d_i -3, \dots, -d_i +1)$, $1 \leq i \leq s$, and we rearrange
it so that a $0$ comes first, followed by the other positive terms in
non-increasing order, followed by the negative terms. Hence, the new
sequence takes the form :
$$
(0,h_1,h_2,\dots,h_r,-h_1,\dots,-h_r)\ ,\ h_1 \geq h_2 \geq \cdots \geq
h_r.
$$
Then, we have :
$$
C(\Omega ) = (h_1 -h_2,h_2-h_3,\dots,h_{r-1}-h_r,h_r).
$$
\item[-] Suppose now that $p+2=2r$. 

If  $(d_1,\dots,d_s)$ does not contain even terms, we use
the same recipe as in the previous case and we obtain the sequence :
$$
(h_1,h_2,\dots,h_r,-h_1,\dots,-h_r)\ , \ h_1 \geq h_2 \geq \cdots \geq
h_r.
$$ 
Then, we have :
$$
C(\Omega) = (h_1-h_2,\dots,h_{r-1}-h_r,h_{r-1}+h_r).
$$

If $(d_1,\dots,d_s)$ contains only even terms. we have the
same sequence :  
$$(
h_1,h_2,\dots,h_r,-h_1,\dots,-h_r)\ ,\ h_1 \geq h_2 \geq \cdots \geq
h_r. 
$$
Then, we set  : 
$a = 0$, if  $r$ is a  multiple of  $4$, and $a=2$, otherwise.
In this situation, there are two orbits $\Omega^{I}$ and $\Omega^{II}$. 
We have :
$$
\begin{array}{c}
C(\Omega^{I}) = (h_1-h_2,\dots, h_{r-2}-h_{r-1},a,2-a),\\
C(\Omega^{II}) = (h_1-h_2,\dots, h_{r-2}-h_{r-1},2-a,a).
\end{array}
$$ 
\end{itemize}

Clearly, a nilpotent element is even if all the terms of its
characteristic are even.

\begin{proof} (of Proposition~\ref{evenness})
Let  $X$ be a nilpotent element of $\mathfrak{p}$, $\Omega_K = K.X$, 
$\Omega = G.X$ and $\Omega_{0}= G_0.X_0$, the corresponding real nilpotent orbit, 
via the Kostant-Sekiguchi correspondence. By~\cite[Remark 9.5.2]{CM}, 
we have $G.X_0 = \Omega$.

Let $(d_1,d_2,\dots,d_s)$ be the Young diagram of $\Omega$. By Remark~\ref{rem1}
and the description of nilpotent orbits in $\mathrm{so}(p,2)$, we only need to
consider the two following cases :
\begin{itemize}
\item[-] $(d_1, \dots,d_s)  \neq (2,2,1,\dots,1)$. Hence, $(d_1,\dots,d_s)$
contains only odd terms. We deduce from the discussions above
that $C(\Omega)$ contains only even terms, and therefore, $X$ is
even.
\item[-] $(d_1, \dots,d_s)  = (2,2,1,\dots,1)$. The two distinct
corresponding real orbits are those of minimal dimension. Defining,
from $(e_{i})$, $1 \leq i \leq 2$, a system of restricted roots for
$\mathfrak{g}_{0}$, we know by a classical result 
(see \cite[\S 4.3]{CM}) that these two orbits
are respectively generated by $X_0 = X_{e_1 - e_2}$ and
$-X_0$. Consider now the element $H_0 = H_{e_1 + e_2}$. Then, $H_0$
is an element of $\mathfrak{p}$ belonging to the centralizer of
$X_0$. Applying the Cayley transform, we check easily that $H_0 \in
\mathfrak{p}^{X}$. So $X$ is not $\mathfrak{p}$-distinguished.
\end{itemize}

We have therefore proved that $\mathfrak{p}$-distinguished elements
are even. 
\end{proof}

\section{Explicit computations}\label{computations}

We shall describe in this section some of the explicit computations
stated in the proofs of Lemmas~\ref{typeB},~\ref{typeD}, and 
Remarks~\ref{smallrank},~\ref{rank2}.
\vskip1em

$\bullet$ Case ($B_{3}$, $\alpha_{1}$). 
We have $S=\{ \alpha_{2},\alpha_{3}\}$, and 
$$
R_{S}^{1}=\{ \alpha_{1} ,\alpha_{1}+\alpha_{2} , \alpha_{1}+\alpha_{2}+\alpha_{3
},
\alpha_{1}+\alpha_{2}+2\alpha_{3},\alpha_{1}+2\alpha_{2}+2\alpha_{3}\}.
$$

A simple computation shows that $X_{\{ \alpha_{1} \} }$
and $X_{\Pi}$ are not $\mathfrak{p}_{S}$-regular. 

For $\lambda \in \mathbb{C}^{*}$, set $X_{\lambda}=X_{\{ \alpha_{1}\}}
+\lambda X_{\Pi}$. Let $Y\in \mathfrak{p}_{S}^{X_{\lambda}}$.

Writing $Y=Y_{1}+Y_{-1}$ with $Y_{1}\in \mathfrak{u}_{S}$ 
and $Y_{-1}\in \mathfrak{u}_{S}^{-}$, we 
conclude that $Y=A+Z$ where $A\in \mathfrak{a}$ and
$Z\in \mathfrak{p}_{S}^{0}=\bigoplus_{\alpha\in R^{0}_{1}} 
\mathfrak{g}_{\alpha}$ with $R_{1}^{0}=(R_{S}^{1}\cup -R_{S}^{1})\setminus
\{ \alpha_{1},-\alpha_{1},\varepsilon_{\Pi}, -\varepsilon_{\Pi }\}$.

So we may assume that $Y\in \mathfrak{p}_{S}^{0}$. 
It follows that
$$
0 = [X_{\lambda},Y] = [X_{\alpha_{1}},Y_{-1}]+\lambda [X_{-\varepsilon_{\Pi}},Y_
{1}]
+ [X_{-\alpha_{1}},Y_{1}] + \lambda [X_{\varepsilon_{\Pi}},Y_{-1}].
$$
So
$$
[X_{\alpha_{1}},Y_{-1}]+\lambda [X_{-\varepsilon_{\Pi}},Y_{1}]
= [X_{-\alpha_{1}},Y_{1}] + \lambda [X_{\varepsilon_{\Pi}},Y_{-1}] = 0.
$$
Hence if $Y=\sum_{\alpha\in R_{1}^{0}} c_{\alpha}X_{\alpha}$,
then $[X_{\lambda},Y]=0$ is equivalent to 
$$
\begin{array}{rl}
c_{-\beta} [X_{\alpha_{1}},X_{-\beta}] + \lambda c_{\gamma} 
[X_{-\varepsilon_{\Pi}},X_{\gamma}] & = 0,\\
c_{\gamma} [X_{-\alpha_{1}},X_{\gamma}] + \lambda c_{-\beta}
[X_{\varepsilon_{\Pi}},X_{-\beta}] & = 0.
\end{array}
$$
whenever we have $\beta,\gamma\in R_{S}^{1}\cap R_{1}^{0}$ verifying
$\beta+\gamma=\alpha_{1}+\varepsilon_{\Pi}$.
The list of such pairs is as follows:
$$
(\beta_{1},\gamma_{1})=(\alpha_{1}+\alpha_{2}, \alpha_{1}+\alpha_{2}+2\alpha_{3}
) \  , \  
(\beta_{2},\beta_{2})=(\alpha_{1}+\alpha_{2}+\alpha_{3}, \alpha_{1}+\alpha_{2}+\alpha_{3}).
$$
Using the Chevalley basis from {\sc Gap4}, we are reduced to the
condition:
$$
\begin{array}{c}
c_{-\beta_{1}}  + \lambda c_{\gamma_{1}} = 0 =
c_{\gamma_{1}}  + \lambda c_{-\beta_{1}},\\
c_{-\beta_{2}}  - \lambda c_{\beta_{2}} = 0 =
c_{\beta_{2}}  - \lambda c_{-\beta_{2}}.
\end{array}
$$
These conditions can only be satisfied if $Y=0$ or $\lambda=\pm 1$.

We have therefore checked that if $X\in \mathfrak{a}\setminus \{ 0\}$ is
not $\mathfrak{p}_{S}$-regular, then it belongs to one of the following
lines : $\mathbb{C} X_{\{ \alpha_{1} \}}$, $\mathbb{C} X_{\Pi}$,
$\mathbb{C} X_{1}$, $\mathbb{C} X_{-1}$.

Using {\sc Gap4}, we obtain that if $X$ belongs to 
the first two lines, then $\dim \mathfrak{g}^{X}=7$,
thus its semisimple part is necessarily of type $A_{1}\times A_{1}$.
If $X$ belongs to the other two lines, then
$\dim \mathfrak{g}^{X}=11$, and its semisimple part must
be of type $B_{2}$.
\vskip1em

$\bullet$ Case ($D_{5},\alpha_{1}$). The elements of $R_{S}^{1}$ are:
$$
\begin{array}{c}
\alpha_{1}, \alpha_{1}+\alpha_{2},\alpha_{1}+\alpha_{2}+\alpha_{3},
\alpha_{1}+\alpha_{2}+\alpha_{3}+\alpha_{4}, \alpha_{1}+\alpha_{2}+\alpha_{3}
+\alpha_{5},\\
\alpha_{1}+\alpha_{2}+\alpha_{3}+\alpha_{4}+\alpha_{5},
\alpha_{1}+\alpha_{2}+2\alpha_{3}+\alpha_{4}+\alpha_{5},\\
\alpha_{1}+2\alpha_{2}+2\alpha_{3}+\alpha_{4}+\alpha_{5}.
\end{array}
$$
Using the same argument above, the pairs are
$$
\begin{array}{c}
(\alpha_{1}+\alpha_{2} , \alpha_{1}+\alpha_{2}+2\alpha_{3}+\alpha_{4}+\alpha_{5}
),\\
(\alpha_{1}+\alpha_{2}+\alpha_{3} ,\alpha_{1}+\alpha_{2}+\alpha_{3}+\alpha_{4}+
\alpha_{5}),\\  
(\alpha_{1}+\alpha_{2}+\alpha_{3}+\alpha_{4},\alpha_{1}+\alpha_{2}+\alpha_{3}+
\alpha_{5}).
\end{array}
$$
Again, we may check using {\sc Gap4} that 
if $X\in \mathfrak{a}\setminus\{ 0\}$ is
not $\mathfrak{p}_{S}$-regular, then it belongs to one of the following
lines : $\mathbb{C} X_{\{ \alpha_{1} \}}$, $\mathbb{C} X_{\Pi}$,
$\mathbb{C} X_{1}$, $\mathbb{C} X_{-1}$.

Using {\sc Gap4}, we obtain that if $X$ belongs to 
the first two lines, then $\dim \mathfrak{g}^{X}=19$,
thus its semisimple part is necessarily of type $A_{1}\times A_{3}$.
If $X$ belongs to the other two lines, then
$\dim \mathfrak{g}^{X}=29$, and its semisimple part must
be of type $D_{4}$.
\vskip1em

$\bullet$ Computations are similar for cases 
$(\mathrm{so}_{6},\mathrm{so}_{4}\times \mathrm{so}_{2})$
and $(\mathrm{so}_{4},\mathrm{so}_{2}\times \mathrm{so}_{2})$.
\vskip1em

$\bullet$ Case ($A_{n},\alpha_{2}$) or ($A_{n},\alpha_{n-1}$) with $n\geq 4$.
Note that the case $n=3$ is just the symmetric pair
$(\mathrm{so}_{6},\mathrm{so}_{4}\times \mathrm{so}_{2})$.
By symmetry, we only need to consider the case $\alpha_{2}$.

Here $\mathcal{E}$ consists of $\Pi$ and
$K=\{ \alpha_{2},\alpha_{3},\dots ,\alpha_{n-1}\}$.
Consider the element $X_{K}$. Its centralizer
contains
$X_{\varepsilon_{\Pi}}$ and $X_{\alpha_{1}+\alpha_{2}}$ because $n\geq 4$.
Since these two root vectors are not
proportional, the subpair $(\mathfrak{r},\mathfrak{r}_{+})$ associated
to $X_{K}$  can not be of type
$(\mathrm{so}_{m+1},\mathrm{so}_{m})$ (see~\cite[Proposition 3]{SY}).

Via {\sc Gap4} computations, we see that the 
semisimple part of $\mathfrak{g}^{X_{K}}$ is of type
$A_{2}$.
\vskip1em

$\bullet$ Case ($D_{5}, \alpha_{4}$) or ($D_{5}, \alpha_{5}$). By symmetry,
we may take $\alpha_{5}$. 
Recall that the two elements of $\mathcal{E}$ are 
$K=\{\alpha_{3},\alpha_{4},\alpha_{5}\}$ and $\Pi$.
Consider the element $X_{K}$. 
Its centralizer contains
$X_{\varepsilon_{\Pi}}$ and $X_{\alpha_{1}+\alpha_{2}+\alpha_{3}+\alpha_{5}}$
which are not proportional. Consequently, 
the subpair $(\mathfrak{r},\mathfrak{r}_{+})$ associated
to $X_{K}$ can not be of type
$(\mathrm{so}_{m+1},\mathrm{so}_{m})$.
\vskip1em

$\bullet$ Case ($E_{6},\alpha_{1}$) or ($E_{6},\alpha_{6}$). Again
we may take $\alpha_{1}$. The set $R_{1}^{+}$ has $16$ elements.
The two elements of $\mathcal{E}$ are
$K=\Pi \setminus \{ \alpha_{2} \}$ and $\Pi$.
Again, consider the element $X_{K}$. 
Its centralizer contains $X_{\varepsilon_{\Pi}}$ and
$X_{\alpha_{1}+\alpha_{2}+\alpha_{3}+\alpha_{4}}$.
So the subpair $(\mathfrak{r},\mathfrak{r}_{+})$ associated
to $X_{K}$ can not be of type
$(\mathrm{so}_{m+1},\mathrm{so}_{m})$.

{\sc Gap4} computations show that the semisimple part of 
$\mathfrak{g}^{X_{K}}$ is of type $A_{5}$. 

\section{A final remark}

In view of Lemma~\ref{typeE7}, we tried to apply the same
method for the symmetric pair $(E_{7},E_{6}\times \mathbb{C})$. 
Unfortunately, using the classification of
$\mathfrak{p}$-distinguished orbits in $\mathfrak{p}$ as described 
in~\cite{PT}, we found a non even $\mathfrak{p}$-distinguished
element. Namely, the orbit corresponding to label 3
of Table 13 of~\cite{PT}.

\begin{acknowledgements}
The authors would like to thank Patrice Tauvel and
Abderrazak Bouaziz for many useful discussions.
\end{acknowledgements}

\affiliationone{Herv\'e Sabourin and
Rupert W.T. Yu \\
UMR 6086 CNRS\\
D\'epartement de Math\'ematiques\\
Universit\'e de Poitiers\\
Boulevard Marie et Pierre Curie\\
T\'el\'eport 2 - BP 30179\\
86962 Futuroscope Chasseneuil cedex\\
France\\
sabourin@math.univ-poitiers.fr\\ 
yuyu@math.univ-poitiers.fr}

\end{document}